\documentclass[12pt]{article}
\usepackage{amsmath,amssymb,amsthm,latexsym,graphicx}
\newtheorem{theorem}{Theorem}

\newtheorem{lemma}[theorem]{Lemma}
\newtheorem{remark}[theorem]{Remark}

\newtheorem{definition}[theorem]{Definition}
\def\r{\rho}
\def\R{\mathcal{R}}
\def\HH{\mathcal{H}}
\def\CC{\Bbb{C}}
\def\RR{\Bbb{R}}
\def\ZZ{\Bbb{Z}}

\begin{document}
\title{A range description for the planar circular Radon transform}
\author{Gaik Ambartsoumian and Peter Kuchment\\
Texas A\& M University, College Station, TX 77843-3368\\
kuchment@math.tamu.edu, haik@tamu.edu}
\date{}
\maketitle
\begin{abstract}
The transform considered in the paper integrates a function
supported in the unit disk on the plane over all circles centered
at the boundary of this disk. Such circular Radon transform arises
in several contemporary imaging techniques, as well as in other
applications. As it is common for transforms of Radon type, its
range has infinite co-dimension in standard function spaces. Range
descriptions for such transforms are known to be very important
for computed tomography, for instance  when dealing with
incomplete data, error correction, and other issues. A complete
range description for the circular Radon transform is obtained.
Range conditions include the recently found set of moment type
conditions, which happens to be incomplete, as well as the rest of
conditions that have less standard form. In order to explain the
procedure better, a similar (non-standard) treatment of the range
conditions is described first for the usual Radon transform on the
plane.
\end{abstract}

\section{Introduction}
The following ``circular'' Radon transform, which is the main
object of study in this article, arises in several applications,
including the newly developing thermoacoustic tomography and its
sibling optoacoustic tomography (e.g.,
\cite{AmbKuch,FPR,Kruger,Patch,MXW,YXW1,YXW2,XWAK}), as well as
radar, sonar, and other applications
\cite{LQ,NC,Norton,Nort-Linzer}. It has also been considered in
relation to some problems of approximation theory, mathematical
physics, and other areas
\cite{ABK,AQ,CH,Leon_Radon,John,Kuch_AMS05,LP1,LP2}.

Let $f(x)$ be a continuous function on $\Bbb{R}^d,\,d\geq 2$.

\begin{definition}\label{D:circular}The circular Radon transform of $f$ is defined
as
$$
Rf(p,\r)=\int_{|y-p|=\r}f(y)d\sigma(y),
$$
where $d\sigma(y)$ is the surface area on the sphere $|y-p|=\r$
centered at $p \in \Bbb{R}^d$.
\end{definition}

In this definition we do not restrict the set of centers $p$ or
radii $r$. It is clear, however, that this mapping is
overdetermined, since the dimension of pairs $(p,r)$ is $d+1$,
while the function $f$ depends on $d$ variables only. This (as
well as the tomographic motivation) suggests to restrict the set
of centers to a set (hypersurface) $S \subset \Bbb{R}^d$, while
not imposing any restrictions on the radii. This restricted
transform is denoted by $R_S$:
$$
R_Sf(p,\r)=Rf(p,\r)|_{p \in S}.
$$

In this paper we will be dealing with the planar case only, i.e.
the dimension $d$ will be equal to $2$. Due to tomographic
applications, where $S$ is the set of locations of transducers
\cite{Kruger,MXW,YXW1,XWAK}, we will be from now on looking at the
specific case when $S$ is the unit circle $|x|=1$ in the plane.

There are many questions one can ask concerning the circular
transform $R_S$: its injectivity, inversion formulas, stability of
inversion, range description, etc. Experience of computerized
tomography shows (e.g., \cite{Natt4,Natt2001}) that all these
questions are of importance. Although none of them has been
resolved completely for $R_S$, significant developments have
occurred recently (e.g.,
\cite{ABK,AQ,AmbKuch,And,Den,Faw,FPR,Gi,KuchQuinto,Natt2001,Nil,Norton,Nort-Linzer,Pal,MXW,YXW1,YXW2}).
The goal of this article is to describe the range of $R_S$ in the
two-dimensional case, with $S$ being the unit circle. Moreover, we
will be dealing with functions supported inside the circle $S$
only. The properties of the operator $R_S$ (e.g., stability of the
inversion, its FIO properties, etc.) deteriorate on functions with
supports extending outside $S$ (e.g., \cite{AQ,FPR,XWAK} and
remarks in the last section of this article). However, in
tomographic applications one normally deals with functions
supported inside $S$ only \cite{Kruger,Patch,MXW,XWAK}.

As it has already been mentioned, the range of $R_S$ has infinite
co-dimension (e.g., in spaces of smooth functions, see details
below) and thus infinitely many range conditions appear. It seems
to be a rather standard situation for various types of Radon
transforms that range conditions split into two types, one of
which is usually easier to discover, while another ``half'' is
harder to come by. For instance, it took about a decade to find
the complete range description for the so called exponential Radon
transform arising in SPECT (single photon emission computed
tomography) \cite{AEK,AK,KL1,KL2,TM}. For a more general
attenuated transform arising in SPECT, it took twice as much time
to move from a partial set of range conditions
\cite{Natt4,Natt-att}  to the complete set \cite{Nov2}. In the
circular case, a partial set of such conditions was discovered
recently \cite{Patch}. It happens to be incomplete, and the goal
of this text is to find the complete one.

One might ask why is it important to know the range conditions.
Such conditions have been used extensively in tomography (as well
as in radiation therapy planning, e.g. \cite{CQ1,CQ2,K1}) for
various purposes: completing incomplete data, detecting and
correcting measurement errors and hardware imperfections,
recovering unknown attenuation, etc.
\cite{Hertle2,Lv,Mennes,Natt3,Natt4,Noo_half,NW,Po2,Sol3,Sol4}.
Thus, as soon as a new Radon type transform arises in an
application, a quest for the range description begins.

In order to explain our approach, we start in the next section
with treating a toy example of the standard Radon transform on the
plane, where the range conditions are well known (e.g.,
\cite{Leon_Radon,GGG1,GGG2,GelfVil,Helg_Radon,Natt4,Natt2001}).
Our approach, however, is different from the standard ones and
naturally leads to the considerations of the circular transform in
the rest of the paper.

\section{The case of the planar Radon transform}

In this section we will approach in a somewhat non-standard way
the issue of the range description for the standard Radon
transform on the plane. Consider a compactly supported smooth
function $f(x)$ on the plane and its Radon transform
\begin{equation}\label{E:Radon}
    (\R f)(\omega,s)=g(\omega,s):=\int\limits_{x\cdot\omega=s}f(x)dl,
\end{equation}
where $s\in\RR$, $\omega\in S^1$ is a unit vector in $\RR^2$, and
$dl$ is the arc length measure on the line $x\cdot\omega=s$. We
want to describe the range of this transform, say on the space
$C^\infty_0 (\RR^2)$. Such a description is well known (e.g.,
\cite{Leon_Radon,GGG1,GGG2,GelfVil,Helg_Radon,Natt4,Natt2001}, or
any other book or survey on Radon transforms or computed
tomography):
\begin{theorem}\label{T:Fourier_range_Radon}
A function $g$ belongs to the range of the Radon transform on
$C^\infty_0$ if and only if the following conditions are
satisfied:
\begin{enumerate}
\item $g\in C_0^\infty(S^1\times \RR)$,

\item for any $k \in \ZZ^+$ the $k$-th moment
$G_k(\omega)=\int\limits_{-\infty}^\infty s^k g(\omega,s)ds$ is
the restriction to the unit circle $S^1$ of a homogeneous
polynomial of $\omega$ of degree $k$,

\item $g(\omega,s)=g(-\omega,-s)$.
\end{enumerate}
\end{theorem}

We would like to look at this result from a little bit different
prospective, which will allow us to do a similar thing in the case
of the circular Radon transform.

In order to do so, let us expand $g(\omega,s)$ into the Fourier
series with respect to the polar angle $\psi$ (i.e., $\omega=(\cos
\psi, \sin \psi)$)
\begin{equation}\label{E:Fourier_series}
g(\omega,s)=\sum\limits_{n=-\infty}^\infty g_n(s)e^{in\psi}.
\end{equation}

We can now reformulate the last theorem in the following a little
bit strange way:
\begin{theorem}\label{T:Fourier_range_Radon_modified}
A function $g$ belongs to the range of the Radon transform on
$C^\infty_0$ if and only if the following conditions are
satisfied:
\begin{enumerate}
\item $g\in C_0^\infty(S^1\times \RR)$,

\item for any $n$, the Mellin transform
$Mg_n(\sigma)=\int\limits_0^\infty s^{\sigma-1} g_n(s)ds$ of the
$n$-th Fourier coefficient $g_n$ of $g$  vanishes at any pole
$\sigma$ of the function $\Gamma(\frac{\sigma+1-|n|}{2})$,

\item $g(\omega,s)=g(-\omega,-s)$.
\end{enumerate}
\end{theorem}

Since the only difference in the statements of these two theorems
is in the conditions 2, let us check that these conditions mean
the same thing in both cases. Indeed, let us expand $g(\omega,s)$
into Fourier series (\ref{E:Fourier_series}) with respect to
$\psi$. Representing  $e^{in\psi}$ as the homogeneous polynomial
$(\omega_1+i(sign\, n) \omega_2)^{|n|}$ of $\omega$ of degree
$|n|$, and noticing that $\omega_1^2+\omega_2^2=1$ on the unit
circle, one easily concludes that the condition 2 in Theorem
\ref{T:Fourier_range_Radon} is equivalent to the following: the
$k$-th moment $\int\limits_{\RR}s^k g_n(s)ds$ of the $n$-th
Fourier coefficient vanishes for integers $0\leq k <|n|$ such that
$k-n$ is even.

Let us now look at the condition 2 in Theorem
\ref{T:Fourier_range_Radon_modified}, still using the same Fourier
expansion. Notice that when $k-|n|$ is a negative even integer,
$Mg_n(\sigma)$ is one-half of the moment of order $k=\sigma-1$ of
$g_n$. Taking into account that
$\Gamma(\frac{\sigma+1-|n|}{2})=\Gamma(\frac{k+2-|n|}{2})$ has
poles exactly when $k-|n|$ is a negative even integer, we see that
conditions 2 in both theorems are in fact saying the same thing.

One can now ask the question, why should one disguise in the
statement of Theorem \ref{T:Fourier_range_Radon_modified} negative
integers as poles of Gamma-function and usual moments as values of
Mellin transforms? The answer is that in the less invariant and
thus more complex situation of the circular Radon transform, one
can formulate a range description in the spirit of Theorem
\ref{T:Fourier_range_Radon_modified}, albeit it is unclear how to
get an analog of the version given in Theorem
\ref{T:Fourier_range_Radon}.

As a warm-up, let is derive the condition 2 in Theorem
\ref{T:Fourier_range_Radon_modified} directly, without relying on
the version given in the preceding theorem. This is in fact an
easy by-product of the A.~Cormack's inversion procedure, see e.g.
\cite[Section II.2]{Natt2001}. Indeed, if we write down the
original function $f(x)$ in polar coordinates $r(\cos \phi,\sin
\phi)$ and expand into the Fourier series with respect to the
polar angle $\phi$
\begin{equation}\label{E:Fourier_series_f}
f(r(\cos \phi,\sin \phi))=\sum\limits_{n=-\infty}^\infty
f_n(r)e^{in\phi},
\end{equation}
then the Fourier coefficients $f_n$ and $g_n$ of the original and
of its Radon transform are related as follows \cite[formula (2.17)
and further]{Natt4}:
\begin{equation}\label{E:Mellin_Radon}
   M(rf_n(r))(s)=\frac{(Mg_n)(s)}{B_n(s)},
\end{equation}
where
\begin{equation}\label{E:B_Gamma}
   B_n(s)=const \frac{\Gamma(s)2^{-s}}{\Gamma((s+1+|n|)/2)\Gamma((s+1-|n|)/2)}
\end{equation}
Thus, condition 2 of Theorem \ref{T:Fourier_range_Radon_modified}
guarantees that the function $M(rf_n(r))(s)$ does not develop
singularities (which it cannot do for a $C_0^\infty$-function $f$)
at zeros of $B_n(s)$. It is not that hard now to prove also
sufficiency in the theorem, applying Cormack's inversion procedure
to $g$ satisfying conditions 1 - 3. However, we are not going to
do so, since in the next sections we will devote ourselves to
doing similar thing in the more complicated situation of the
circular Radon transform.

\section{The circular Radon transform. Formulation of the main result}

Let us recall the notion of Hankel transform (e.g.,
\cite{Davies}). For a function $h(r)$ on $\RR^+$, one defines its
Hankel transform of an integer order $n$ as follows:
\begin{equation}\label{E:Hankel_transform}
(\HH_n h)(\sigma)=\int\limits_0^\infty J_n(\sigma r)h(r)r\,dr,
\end{equation}
where the standard notation $J_n$ is used for Bessel functions of
the first kind.

Let, as in the Introduction, $R_S$ be the circular Radon transform
on the plane that integrates functions compactly supported inside
the unit disk $D$ over all circles $|x-p|=\r$ with centers $p$
located on the unit circle $S=\{p\,|\,|p|=1$). Since this
transform commutes with rotations about the origin, the Fourier
series expansion with respect to the polar angle partially
diagonalizes the operator, and thus the $n$-th Fourier coefficient
$g_n(\r)$ of $g=R_S f$ will depend on the $n$-th coefficient $f_n$
of the original $f$ only. It was shown in \cite{Norton} that the
following relation between these coefficients holds:
\begin{equation}\label{E:Norton}
  g_n(\r)=2\pi \r \HH_0\{J_n\HH_n\{f_n\}\}.
\end{equation}
For the reader's convenience, we will provide the brief derivation
from \cite{Norton}. Considering a single harmonic
$f=f_n(r)e^{in\phi}$ and using polar coordinates, one obtains
\begin{equation}\label{E:intermed}
g_n(\r)=\int\limits_0^\infty rf_n(r)dr\int\limits_0^{2\pi}
\delta\left[(r^2+1-2r\cos \phi)^{1/2}-\rho\right]e^{-in\phi}d\phi.
\end{equation}
Thus, the computation boils down to evaluating the integral
$$
I=\int\limits_0^{2\pi} \delta\left[(r^2+1-2r\cos
\phi)^{1/2}-\rho\right]e^{-in\phi}d\phi.
$$
Using the standard identity
$$
\delta(\rho^\prime-\rho)=\rho\int\limits_0^\infty J_0(\rho^\prime
z)J_0(\rho z)zdz
$$
and the identity that is easy to obtain from one of the addition
formulas, e.g. from \cite[formula (4.10.6)]{Andrews}
$$
2\pi J_n(a z)J_n(b z)=\int\limits_0^{2\pi} J_0[z(a^2+b^2-2ab\cos
\phi)^{1/2}]e^{-in\phi}d \phi,
$$
one arrives from (\ref{E:intermed}) to (\ref{E:Norton}).

Since Hankel transforms are involutive, it is easy to invert
(\ref{E:Norton}) and get Norton's inversion formulae \cite{Norton}
\begin{equation}\label{E:Norton_inversion}
    f_n=\frac{1}{2\pi}\HH_n\{\frac{\HH_0\{g_n(\r)/\r\}}{J_n}\}.
\end{equation}
Now one can clearly see analogies with the case of the Radon
transform, where zeros of Bessel functions should probably
introduce some range conditions. This happens to be correct and
leads to the main result of this article:
\begin{theorem}\label{T:main}
In order for the function $g(p,\r)$ on $S^1\times\RR$ to be
representable as $R_S f$ with $f\in C_0^\infty(D)$, it is necessary
and sufficient that the following conditions are satisfied:
\begin{enumerate}
\item $g\in C_0^\infty(S^1\times (0,2))$.

\item For any $n$, the $2k$-th moment $\int\limits_0^\infty
\r^{2k} g_n(\r)d\r$ of the $n$-th Fourier coefficient of $g$
vanishes for integers $0\leq k <|n|$. (Equivalently, the $2k$-th
moment $\int\limits_0^\infty \r^{2k} g(p,\r)d\r$ is the
restriction to the unit circle $S$ of a (non-homogeneous)
polynomial of $p$ of degree at most $k$.)

\item For any $n\in\ZZ$, function $\HH_0
\{g_n(\r)/\r\}(\sigma)=\int\limits_0^\infty J_0(\sigma
\r)g_n(\r)d\r$ vanishes at any zero $\sigma\neq 0$ of Bessel
function $J_n$. (Equivalently, the $n$th Fourier coefficient with
respect to $p\in S^1$ of the ``Bessel moment'' $G_\sigma
(p)=\int\limits_0^\infty J_0(\sigma \r)g(p,\r)d\r$ vanishes if
$\sigma\neq 0$ is a zero of Bessel function $J_n$.)
\end{enumerate}
\end{theorem}

\section{Proof of the main result}
Let us start with proving necessity, which is rather
straightforward. Indeed, the necessity of condition 1 is obvious.
Let us prove the second condition. In fact, it has already been
established in \cite{Patch}. Let us repeat for completeness its
simple proof. Let $k$ be an integer. Consider the moment of order
$2k$ of $g$:
\begin{equation}\label{E:necess_moment}
    \int\limits_0^\infty \r^{2k}
    g(p,\r)d\r=\int\limits_{\RR^2}|x-p|^{2k}f(x)dx=\int\limits_{\RR^2}(|x|^2-2x\cdot p
    +1)^{k}f(x)dx
\end{equation}
(we have taken into account that $|p|=1$). We see that the
resulting expression is the restriction to $S^1$ of a
(non-homogeneous) polynomial of degree $k$ in variable $p$.
Expanding into Fourier series with respect to the polar angle of
$p$, we see that the $n$th harmonic $g_n$ contributes the
following homogeneous polynomial of degree $|n|$ in the variable
$p$:
$$
\left(\int\limits_0^\infty \r^{2k}g_n(\r)d\r\right)e^{in\psi}.
$$
Here as before $p=(\cos \psi, \sin \psi)$. Thus, for $|n|>k$, this
term must vanish, which gives necessity of condition 2. We will
return to a discussion of this condition below to add a new twist
to it.

Necessity of condition 3 follows immediately from Norton's formula
(\ref{E:Norton_inversion}), which implies in particular that
$$
\HH_0\{g_n(\r)/\r\}=2\pi J_n\HH_n\{f_n\}.
$$
Since both functions $J_n$ and $\HH_n\{f_n\}$ are entire,
$\HH_0\{g_n(\r)/\r\}$ vanishes whenever $J_n$ does.

\begin{remark}\label{R:zeros} The reader might ask why in the third condition of the
Theorem we do not take into account the zero root of $J_n$, which
in fact has order $n$, while non-zero roots are all simple. The
reason is that the condition 2 already guarantees that $\sigma=0$
is zero of order $2n$ of $\HH_0\{g_n(\r)/\r\}$ (twice higher than
that of $J_n$). Indeed, due to evenness of $J_0$, function
$\HH_0\{g_n(\r)/\r\}(\sigma)$ is also even. Thus, all odd order
derivatives at $\sigma=0$ vanish. The known Taylor expansion of
$J_0$ at zero leads to the formula
$$
\HH_0\{g_n(\r)/\r\}(\sigma)=\sum\limits_m\frac{(-1)^m}{(m!)^2}
\left(\frac{\sigma}{2}\right)^{2m}\int\limits_0^\infty
r^{2m}g_n(r)dr.
$$
We see now that the moment condition 2 guarantees that $\sigma=0$
is zero of order $2n$ of $\HH_0\{g_n/\r\}(\sigma)$.
\end{remark}

Let us move to the harder part, proving sufficiency. Assume a
function $g$ satisfies conditions of the theorem and is supported
in $S\times (\epsilon,2-\epsilon)$ for some positive $\epsilon$.
We will show that then $g=R_Sf$ for some $f\in C^\infty_0
(D_\epsilon)$, where $D_\epsilon$ is the disk $|x|<\epsilon$ in
the plane.

Due to Norton's formulas, it is natural to expect the proof to go
along the following lines: expand $g$ into the Fourier series
$g=\sum\limits_m g_m(\r)e^{im\psi}$ with respect to the angle
variable $\psi$, then use (\ref{E:Norton_inversion}) to construct
a function $f$ and then show that $f$ is of an appropriate
function class and that its circular Radon transform is equal to
$g$. This is what we are going to do, with a small caveat that
instead of constructing $f$ itself, we will construct its
two-dimensional Fourier transform. Besides, we will start
considering the partial sums of the series
$h_n=\sum\limits_{|m|\leq n} g_m(\r)e^{im\phi}$. But first, we
need to get some simple estimates from below for the Bessel
function of the first kind $J_n$.

\begin{lemma}\label{L:bessel}
On the entire complex plane except for a disk $S_0$ centered at
the origin and a countable number of disks $S_k$ of radii $\pi/6$
centered at points $\pi(k+\frac{2n+3}{4})$, one has

\begin{equation}
\label{E:estimate} |J_n(z)|\ge\frac{C
e^{|Im\,z|}}{\sqrt{|z|}},\;\;\;\;\; C>0
\end{equation}

\end {lemma}

\textbf{Proof:} Let us split the complex plane into three parts by
a circle $S_0$ of a radius $R$ (to be chosen later) centered at
the origin and a planar strip $\{z=x+iy|\,|y|<a\}$, as follows:
part I consists of points $z$ satisfying $|z|\geq R$ and
$|Im\,z|\geq a$; part II consists of points such that $|z|\geq R$
and $|Im\,z|< a$; part III is the interior of $S_0$, i.e. $|z|<R$.
It is clearly sufficient to prove the estimate (\ref{E:estimate})
in the first two parts: outside and inside the strip. Using the
parity property of $J_n$, it suffices to consider only the right
half plane $Re\,z\geq 0$.

The Bessel function of the first kind $J_n(z)$ has the following
known asymptotic representation in the sector $|\arg
z|\le\pi-\delta$ (e.g., \cite[formula (4.8.5)]{Andrews} or
\cite[formula (5.11.6)]{Lebedev}):
\begin{equation}\label{E:asymptotic}
\begin{array}{c}
J_n(z)=\sqrt{\frac{2}{\pi z}} \cos(z-\frac{\pi n}{2}
-\frac{\pi}{4})(1+O(|z|^{-2}))\\
 -\sqrt{\frac{2}{\pi z}} \sin(z-\frac{\pi
n}{2}-\frac{\pi}{4})\left( \frac{4n^2-1}{8z}+O(|z|^{-3})\right)
\end{array}
\end{equation}

Let us start estimating in the first part of the complex plane,
i.e. where $|Im \, z|>a$ and $|z|>R$ for sufficiently large $a$
and $R$ (and, as we have agreed, $Re\,z \geq 0$). There, due to
boundedness of $\tan z$ in this region, one concludes that
$\displaystyle{\frac{\sin z}{z} = \cos z\; ( O(|z|^{-1}))}$, and
thus (\ref{E:asymptotic}) implies
$$
J_n(z)=\sqrt{\frac{2}{\pi z}} \cos(z-\frac{\pi n}{2}
-\frac{\pi}{4})(1+O(|z|^{-1})),
$$
which in turn for sufficiently large $a, R$ leads to
\begin{equation}
\label{E:part1} |J_n(z)|\ge\frac{C e^{|Im\, z|}}{\sqrt {|z|}}
\end{equation}

In the second part of the plane (right half of the strip), due to
boundedness of $\sin z$ we have
$$
J_n(z)=\sqrt{\frac{2}{\pi z}}\left[\cos(z-\frac{\pi n}{2}
-\frac{\pi}{4})(1+O(|z|^{-2}))+ O(|z|^{-1})\right].
$$
Consider the system of non-intersecting circles $S_k$ with centers
at $z_k=\frac{\pi}{2}+k\pi+\frac{\pi n}{2}+\frac{\pi}{4}$ and
radii equal to $\frac{\pi}{6}$. Then outside these circles
$|\cos(z-\frac{\pi n}{2} -\frac{\pi}{4})|\ge C$ and
$$
|J_n(z)|\ge
\frac{C}{\sqrt{|z|}}(1+O(|z|^{-1})).
$$
This implies that for a suitably chosen and sufficiently large
$R$, inside of the strip and outside the circles $S_k$, we have
\begin{equation}
\label{E:part2} |J_n(z)|\ge\frac{C e^{|Im\, z|}}{\sqrt{|z|}}
\end{equation}
for $|z|>R$. This proves the statement of the lemma. \qed

Let us now return to our task: consider the function $g$ and the
partial sums $h_n$ of its Fourier series.

\begin{lemma}\label{L:suffic finite}
\begin{enumerate}
\item If $g(\phi,\r)=\sum_m g_m(\r)e^{im\psi}$ satisfies
conditions of Theorem \ref{T:main} and is supported in $S\times
(\epsilon,2-\epsilon)$, then each partial sum $h_n=\sum_{|m|<n}
g_m(\r)e^{im\psi}$ does so.

\item For any $n$, $h_n$ is representable as $R_S f_n$ for a
function $f_n\in C^\infty_0 (D_\epsilon)$.
\end{enumerate}
\end{lemma}
{\bf Proof of the lemma.} The first statement of the lemma is obvious.

Due to i), it is sufficient to prove the second statement for a
single term $g=g_n(\r)e^{in\psi}$. As it was just mentioned, we
will reconstruct the Fourier transform $F$ of the function $f$. In
order to do this, we will use the standard relation between
Fourier and Hankel transforms. Let as before
$f(x)=f_n(r)e^{in\phi}$, where $r=|x|$ and $\phi$ are polar
coordinates on $\RR^2$. Then the Fourier transform $F(\xi)$ of $f$
at points of the form $\xi=\sigma \omega$, where $\sigma\in\CC$
and $\omega=(\cos \psi, \sin \psi)\in\RR^2$ can be written up to a
constant factor as follows:
\begin{equation}\label{E:Fourier-Hankel}
    F(\sigma \omega)=\HH_n(f_n)(\sigma) e^{in\psi}.
\end{equation}
(e.g., \cite[end of Section 14.1]{Davies}). If we knew that $g=R_S
f$, then according to (\ref{E:Norton}) this would mean that
\begin{equation}\label{E:F_defined}
F(\sigma
\omega)=F(\sigma)e^{in\psi}=\frac{1}{2\pi}\frac{\HH_0(g_n(\r)/\r)(\sigma)}{J_n(\sigma)}e^{in\psi}.
\end{equation}
Let us now take this formula (\ref{E:F_defined}) as the definition
of $F(\sigma \omega)$. Due to the standard parity property of
Bessel functions, such $F$ is a correctly defined function of
$\sigma \omega$ for $\sigma\neq 0$ (i.e.,
$F(\sigma\omega)=F((-\sigma)(-\omega))$). We would like to show
that it is the Fourier transform of a function $f\in C^\infty_0
(D_\epsilon)$. Let us prove first that $F$ belongs to the Schwartz
space $\mathcal{S}(\RR^2)$. In order to do so, we need to show its
smoothness with respect to the angular variable $\psi$, smoothness
and fast decay with all derivatives in the radial variable
$\sigma$, as well as that no singularity arises at the origin,
which in principle could, due to usage of polar coordinates.
Smoothness with respect to the angular variable is obvious, due to
(\ref{E:F_defined}). Let us deal with the more complex issue of
smoothness and decay with respect to $\sigma$. First of all,
taking into account that $g_n(\r)$ is supported inside $(0,2)$,
and due to the standard $2D$ Paley-Wiener theorem, we conclude
that $u(\sigma)=\HH_0(g_n(\r)/\r)$ is an entire function that
satisfies for any $N$ the estimate
\begin{equation}\label{E:PWestimate}
|u(\sigma)|\leq C_N (1+|\sigma|)^{-N}e^{(2-\epsilon)|Im\,\sigma|}.
\end{equation}
According to the range conditions 2 and 3 of the Theorem, this
function vanishes at all zeros of Bessel function $J_n (\sigma)$
at least to the order of the corresponding zero. This means, that
function $F(\sigma)=\dfrac{u(\sigma)}{2\pi J_n(\sigma)}$ is
entire. Let us show that it belongs to a Paley-Wiener class.
%

Indeed, $\HH(g_n(\r)/\r)$ is an entire function with Paley-Wiener
estimate (\ref{E:PWestimate}). Due to the estimate from below for
$J_n$ (\ref{E:estimate}) given in Lemma \ref{L:bessel}, we
conclude that $F(\sigma\omega)$ is an entire function of
Paley-Wiener class in the radial directions, uniformly with
respect to the polar angle. Namely,

\begin{equation}\label{E:F_PWestimate}
|F(\sigma)|\leq C_N (1+|\sigma|)^{-N}e^{(1-\epsilon)|Im\,\sigma|}.
\end{equation}

Indeed, outside the family of circles $S_k$ the estimate
(\ref{E:estimate}) together with (\ref{E:PWestimate}) give the
Paley-Wiener estimate (\ref{E:F_PWestimate}) we need. Inside these
circles, application of the maximum principle finishes the job.
Smoothness with respect to the polar angle is obvious. Thus, the
only thing one needs to establish to verify that $F$ belongs to
the Schwartz class is that $F$ is smooth at the origin. This,
however, is the standard question in the Radon transform theory,
the answer to which is well known (e.g., \cite[p. 108--109]{GGG1},
\cite{GGG2,GelfVil}, \cite[Ch. 1, proof of Theorem
2.4]{Helg_Radon}). Namely, one needs to establish that for any
non-negative integer $k$, the $k$th radial (i.e., with respect to
$\sigma$) derivative of $F(\sigma\omega)$ at the origin is a
homogeneous polynomial of order $k$ with respect to $\omega$. So,
let us check that this condition is satisfied in our situation.
First of all, the parity of the function $F$ is the same as of
$n$. Thus, we do not need to worry about the derivatives
$F^{(k)}_\sigma|_{\sigma=0}$ with $k-n$ odd, since they are zero
automatically. Due to the special single-harmonic form of $F$, we
only need to check that $F^{(k)}_\sigma|_{\sigma=0}=0$ for $k<|n|$
with $k-n$ even.This, however, as we have discussed already in
Remark \ref{R:zeros}, follows from the moment conditions 2 of the
Theorem.

Due to the smoothness that we have just established and
Paley-Wiener estimates, $F\in \mathcal{S}(\RR^2)$. Thus,
$F=\hat{f}$ for some $f\in \mathcal{S}(\RR^2)$. It remains to show
that $f$ is supported inside the disk $D_\epsilon$. Consider the
usual Radon transform $\R f(s,\phi)$ of $f$. According to the
standard Fourier-slice theorem
\cite{Leon_Radon,GGG1,GGG2,GelfVil,Helg_Radon,Natt4}, the
one-dimensional Fourier transform (denoted by a ``hat'') from the
variable $s$ to $\sigma$ gives (up to a fixed constant factor) the
values $\widehat{\R f}(\sigma,\psi)=F(\sigma\omega)$, if as before
$\omega=(\cos \psi, \sin \psi)$. Here $\R$, as before, denotes the
standard Radon transform in the plane. Since functions
$F(\sigma\omega)$ of $\sigma$, as we have just discussed, are
uniformly with respect to $\omega$ of a Paley-Wiener class, this
implies that $\R f(s,\omega)$ has uniformly with respect to
$\omega$ bounded support in $|s|<1-\epsilon$. Now the ``hole
theorem'' \cite{Helg_Radon,Natt4} (which is applicable to
functions of the Schwartz class), implies that $f$ is supported in
$D_\epsilon$.

The last step is to show that $R_S f=g=g_n(r)e^{in\phi}$. This,
however, immediately follows from comparing formulas
(\ref{E:F_defined}) and (\ref{E:Norton}), which finishes the proof
of the main Lemma \ref{L:suffic finite}\qed

Let us now return to the proof of Theorem \ref{T:main}. We have
proven so far that any partial sum $h_n$ of the Fourier series for
$g$ belongs to the range of the operator $R_S$ acting on smooth
functions supported inside the disk $D_\epsilon$. The function $g$
itself is the limit of $h_n$ in $C^\infty_0(S\times
(\epsilon,2-\epsilon))$. The only thing that remains to be proven
is that the range is closed in an appropriate topology. Microlocal
analysis can help with this.

Consider $R_S$ as an operator acting from functions defined on the
open unit disk $D$ to functions defined on the open cylinder
$\Omega=S\times (0,2)$. As such, it is a Fourier integral operator
\cite{Guillemin,Guill_Ster,Q1980}. If $R^t_S$ is the dual
operator, then $E=R^t_S R_S$ is an elliptic pseudo-differential
operator of order $-1$ \cite[Theorem
1]{Guillemin}\footnote{Bolker's injective immersion condition
\cite{Guillemin} that is needed for validity of this result, is
satisfied here, as shown in the proof of Lemma 4.3 in \cite{AQ}.}.

\begin{lemma}\label{L:semi-Fredholm}
The continuous linear operator $E:H^2_0 (D_\epsilon)\mapsto
H^3_{loc}(D)$ has zero kernel and closed image.
\end{lemma}
{\bf Proof of the lemma.} Since $E=R^t_S R_S$, the kernel of this
operator coincides with the kernel of $R_S$ acting on $H^2_0
(D_\epsilon)$. Since $S$ is closed, it is known that $R_S$ has no
compactly supported functions in its kernel \cite{ABK,AQ} (this
also follows from analytic ellipticity of $E$ and Theorem 8.5.6 of
\cite{Hormander}, see also Lemma 4.4 in \cite{AQ}). Thus, the
statement about the kernel is proven and we only need to prove the
closedness of the range.

Let $P$ be a properly supported pseudo-differential parametrix of
order $1$ for $E$ \cite{Shubin}. Then $PE=I+B$, where $B$ is an
infinitely smoothing operator on $D$. Consider the operator $\Pi$
that acts as the composition of restriction to $D_\epsilon$ and
then orthogonal projection onto $H^2_0(D_\epsilon)$ in
$H^2(D_\epsilon)$. On $H^2_0(D_\epsilon)$ one has $\Pi P E=I+K$,
where $K$ is a compact operator on $H^2_0(D_\epsilon)$. Notice
that the operator $\Pi P$ is continuous from the Frechet space
$H^3_{loc}(D)$ to $H^2_0(D_\epsilon)$. Due to the Fredholm
structure of the operator $\Pi P E=I+K$ acting on
$H^2_0(D_\epsilon)$, its kernel is finite-dimensional. Let
$M\subset H^2_0(D_\epsilon)$ be a closed subspace of finite
codimension complementary to the kernel, so $I+K$ is injective on
$M$ and has closed range. Then one can find a bounded operator $A$
in $H^2_0(D_\epsilon)$ such that $A(I+K)$ acts as identity on $M$.
Thus, the operator $A\Pi P$ provides a continuous left inverse to
$E:M \mapsto H^3_{loc}(D)$. This shows that the range of $E$ on
$M$ is closed in $H^3_{loc}(D)$. On the other hand, the total
range of $E$ differs only by a finite dimension from the one on
$M$. Thus, it is also closed. \qed

We can now finish the proof of the theorem. Indeed, the last lemma
shows that the function $R^t_S g$, being in the closure of the
range, is in fact in the range, and thus can be represented as
$Ef$ with some $f\in H^2_0(D_\epsilon)$. In other words,
$R^t_S(R_S f - g)=0$. Since the kernel of $R^t_S$ on compactly
supported functions is orthogonal to the range of $R_S$, we
conclude that $R_S f - g=0$. Since $Ef=R^t_S g$ is smooth, due to
ellipticity of $E$ we conclude that $f$ is smooth as well. This
concludes the proof of the theorem. \qed

\section{Remarks and acknowledgments}
We would like to finish with some remarks.

\begin{itemize}
\item It should be possible to prove that the operator $R_S$ in
the situation considered in the text is semi-Fredholm between
appropriate Sobolev spaces (analogously to the properties of the
standard and attenuated Radon transforms, e.g.
\cite{Heike,Natt4}). This would eliminate the necessity of the
closedness of the range discussion in the end of the proof of
Theorem \ref{T:main}.

Such a statement could probably be proven either by using FIO
techniques, or by controlling dependence on $n$ of the constant
$C$ and of the radius of the circle $S_0$ in Lemma \ref{L:bessel}.
The former approach would be better, being more general.

\item Proving compactness of support of function $f$ in Lemma
\ref{L:suffic finite}, we used the standard Radon transform and
the ``hole theorem.'' Instead, one could probably use the fact
that Fourier transform of $f$ is, by construction, a Paley-Wiener
class CR-function on the three-dimensional variety of points
$\sigma\omega$ in $\CC^2$ and then use an appropriate mandatory
analytic extension theorem in the spirit of \cite{Oktem}.

\item We considered the situation most natural for tomographic
imaging, when the functions to reconstruct are supported inside
the aperture curve $S$. What happens when the supports of
functions extend outside the circle $S$? It is known that
compactly supported \cite{AQ} (or even belonging to $L_p$ with
sufficiently small $p$ \cite{ABK}) functions can still be uniquely
reconstructed. Necessity of the range conditions we derived
apparently still holds and they are still sufficient for finite
Fourier series. However, many things do go wrong in this case. Our
proof of the closedness of the range fails (in particular, since
the Bolker's condition for the corresponding FIO does not hold
anymore, which was also the main hurdle in proving the results of
\cite{AQ}). Moreover, the range will not be closed anymore.
Indeed, reconstruction will become unstable, since due to standard
microlocal reasons \cite{KLM,KuchQuinto,LQ,Q1993,XWAK}, some parts
of the wave front set of the function outside $S$ will not be
stably recoverable. This means, in particular, that non-smooth
functions can have smooth circular Radon images. This, in turn
implies that the range is not closed in the spaces under
consideration, and so sufficiency of the range conditions should
fail.  We are not sure what kind of range description, if any,
could work in this situation. By the way, the nice backprojection
type inversion formulas available in odd dimensions \cite{FPR}
also fail for such functions.

\item It would be interesting to understand range conditions in
the case of a closed curve $S$ different from a circle. Since our
method uses rotational invariance, it is not directly applicable
to this situation.

\item Our result is stated and proven in $2D$ only. It is possible
that a similar approach might work in higher dimensions. As we
have been notified by D.~Finch, he and Rakesh have recently
obtained by different methods some range descriptions in $3D$
\cite{FinRak}.

\item The result of the paper was presented at the Fully
Three-Dimensional Image Reconstruction Meeting in Radiology and
Nuclear Medicine July 6-9, 2005 in Salt Lake City, Utah.

\end{itemize}

This work was supported in part by the NSF Grants DMS 9971674 and
0002195. The authors thank the NSF for this support. Any opinions,
findings, and conclusions or recommendations expressed in this
paper are those of the authors and do not necessarily reflect the
views of the National Science Foundation.

\end{document}